\documentclass[12pt,reqno]{article}

\usepackage[usenames]{color}
\usepackage{graphicx}
\usepackage{amscd}
\usepackage{color}
\usepackage{fullpage}
\usepackage{psfig}
\usepackage{graphics,amsmath,amssymb}
\usepackage{amsthm}
\usepackage{amsfonts}
\usepackage{latexsym}
\usepackage{epsf}
\usepackage{float}
\usepackage{bm}

\usepackage[colorlinks=true,
linkcolor=webgreen,
filecolor=webbrown,
citecolor=webgreen]{hyperref}

\definecolor{webgreen}{rgb}{0,.5,0}
\definecolor{webbrown}{rgb}{.6,0,0}

\newcommand{\seqnum}[1]{\href{http://oeis.org/#1}{\underline{#1}}}
\newcommand{\beql}[1]{\begin{equation}\label{#1}}
\newcommand{\eeq}{\end{equation}}
\newcommand{\eqn}[1]{(\ref{#1})}

\newcommand{\NN}{\mathbb N}
\newcommand{\ZZ}{\mathbb Z}

\newtheorem{thm}{Theorem}{\bfseries}{\itshape}
{\bfseries}{\itshape}
{\bfseries}{\itshape}
{\bfseries}{\itshape}
\newtheorem{conj}[thm]{Conjecture}{\bfseries}{\itshape}

\newcommand{\sA}{{\mathcal A}}

\begin{document}
\theoremstyle{plain}

\begin{center}
{\large\bf The On-Line Encyclopedia of Integer Sequences } \\
\vspace*{+.2in}

N. J. A. Sloane\footnote{Neil J. A. Sloane is 
president of the OEIS Foundation and a visiting scientist at Rutgers Univerity.
His email address is njasloane@gmail.com. } 
\ \\
May 25, 2018 \\
\ \\
(To appear in {\em Notices of American Mathematical Society})

\vspace{.3 in}

{\bf Abstract}

\end{center}

The recent history of {\em The On-Line Encyclopedia of Integer Sequences} (or {\em OEIS}), describing developments since 2009, and discussing recent sequences involving interesting unsolved problems and in many cases spectacular illustrations. These include: Peaceable Queens, circles in the plane, the earliest cube-free binary sequence, the EKG and Yellowstone permutations, other lexicographically earliest sequences, iteration of number-theoretic functions, home primes and power trains, a memorable prime, a missing prime, Post's tag system, and coordination sequences.

\section{Introduction}\label{Sec1}

The {\em OEIS}  (or {\em On-Line Encyclopedia of 
Integer Sequences}\footnote{\url{http://oeis.org}.})
is a freely accessible database of number sequences,  now
in its $54$th year, and  online since $1995$.
It contains over $300,000$ entries,
and for each one gives a definition, properties, references, computer programs,
tables, etc., as appropriate. It is widely referenced: a web
page\footnote{\url{http://oeis.org/wiki/Works_Citing_OEIS}.} lists over $6000$
works that cite it, and often say  things like
``this theorem would not exist without the help of the OEIS''.
It has been called one of the most useful mathematical sites in the web.

The main use is to serve as a dictionary or fingerprint file for identifying
number sequences (and when you find the sequence you
are looking for, you will understand why the OEIS is so popular).
If your sequence is not recognized, you see a message saying that if 
the sequence is of general interest, you should
submit it for inclusion in the database. The resulting queue of new submissions is
a continual source of lovely problems.

I described the OEIS in a short article in the September $2003$ issue of these {\em Notices}.
The most significant changes since then took place in $2009$, when
a non-profit foundation\footnote{The OEIS Foundation, Inc.,
\url{http://oeisf.org}.} was set up to own and maintain the OEIS,
and in $2010$ when  the OEIS was moved off my home page at AT\&T Labs to
a commercial host. The format has also changed: since 2010 the OEIS has been
a refereed ``wiki''.
Four people played a crucial role in the transition: Harvey P. Dale and Nancy C. Eberhardt 
helped set up the Foundation, Russell S. Cox wrote the software, 
and David L. Applegate  helped move the OEIS.  The OEIS would probably
not exist today but for their help.

All submissions, of new sequences and updates, are now refereed 
by volunteer editors. One of the rewards of being an editor
is that you see a constant flow of new problems, often 
submitted by non-mathematicians, which frequently contain juicy-looking questions
that are begging to be investigated.

This article will describe a selection of recent sequences, mostly connected 
with unsolved problems.

Sequences in the OEIS are identified by a $6$-digit number
prefixed by~A. \seqnum{A000001} is the number of
groups of order $n$, \seqnum{A000002} is Kolakoski's sequence, and so on.
When we were approaching a quarter of a million entries,
the editors voted to decide which sequence would become \seqnum{A250000}.
The winner was the Peaceable Queens sequence, described in the next section,
and the runner-up was the ``circles in the plane'' sequence \seqnum{A250001}
 discussed in~\S\ref{SecWild}.
The $n$th term of the sequence under discussion is usually denoted by $a(n)$.


\begin{figure}[!ht]
\centerline{\includegraphics[angle=0, width=3.0in]{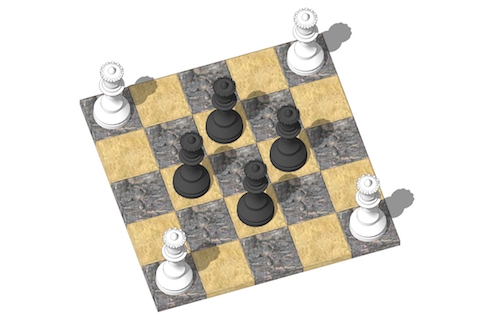}}
\caption{One of three solutions to the Peaceable Queens problem on a $5 \times 5$ board,
illustrating $a(5)=4$.}
\label{FigPQ5}
\end{figure}

\begin{figure}[!ht]
\centerline{\includegraphics[angle=0, width=3.0in]{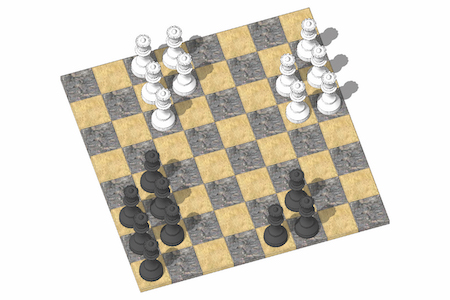}}
\caption{A solution to the Peaceable Queens problem on an $8 \times 8$ board, illustrating $a(8)=9$.
(There are actually $10$ white queens here
but only $9$ count since the numbers of white and black queens must be equal. Any one of the white queens could be omitted.)}
\label{FigPQ8}
\end{figure}

\section{Peaceable Queens}\label{SecPQ}

\begin{figure}[!ht]
\centerline{\includegraphics[angle=0, width=5.0in]{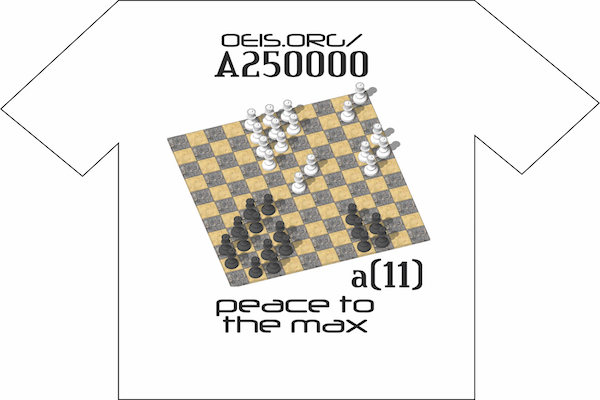}}
\caption{A solution to the Peaceable Queens problem on an $11 \times 11$ board, illustrating $a(11)=17$.}
\label{FigPQ11}
\end{figure}

\begin{figure}[!ht]
\centerline{\includegraphics[angle=0, width=3.0in]{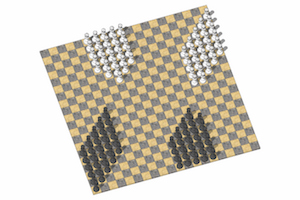}}
\caption{A conjectured solution to the Peaceable Queens problem on a $20 \times 20$ board,
found by Bob Selcoe,
showing that $a(20) \ge 58$.}
\label{FigPQ20}
\end{figure}

In \seqnum{A250000}, $a(n)$
 is the maximal number $m$ such that it is possible
to place $m$ white queens and $m$ black queens on an $n \times n$ chess board so that no 
queen attacks a queen of the opposite color.
These are peaceable queens.
This is a fairly new problem with
some striking pictures, an interesting conjecture,
and a satisfactorily non-violent theme.
It was posed by Robert A. Bosch in 1999,
as a variation on the classical problem
of finding the number of ways to place $n$ queens 
on an $n \times n$ board so that they do not attack each other (\seqnum{A000170}).
It was added to the OEIS in 2014 by Donald E. Knuth,
and a number of people have contributed to the entry since then.
Only thirteen terms are known:
$$
\begin{array}{crrrrrrrrrrrrrr}
n: & 1 & 2 & 3 & 4 & 5 & 6 & 7 & 8 & 9 & 10 & 11 & 12 & 13 \\
a(n): & 0 & 0 & 1 & 2 & 4 & 5 & 7 & 9 & 12 & 14 & 17 & 21 & 24
\end{array}
$$ 
Figures\footnote{The illustrations here are of poor quality because of the limits on file size imposed by the {\em arXiv}. For publication-quality illustrations see the version of this article on the author's home page \url{http://NeilSloane.com}.}
\ref{FigPQ5}-\ref{FigPQ20} show examples of
solutions for $n = 5, 8, 11$ and (conjecturally) $20$.

\begin{figure}[!ht]
\centerline{\includegraphics[angle=0, width=3.0in]{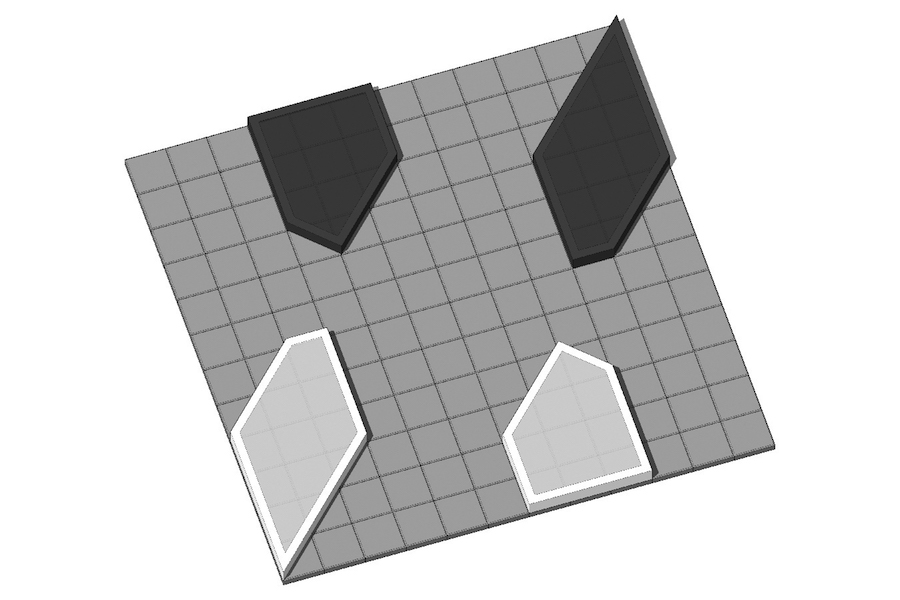}}
\caption{A general construction for the Peaceable Queens problem 
found by Beno\^{\i}t~Jubin, showing that for large $n$, 
$a(n) \ge \lfloor 7 n^2/48 \rfloor$,
a formula which might be exact for all $n > 9$.}
\label{FigPQPK}
\end{figure}

For larger values of $n$, the best solutions presently known were found by Beno\^{\i}t~Jubin
and concentrate the queens into four pentagonal regions,
as shown in Figure~\ref{FigPQPK} (and
generalize the arrangement shown in Figure~\ref{FigPQ20}).
This construction gives a lower bound of $\lfloor 7 n^2/48 \rfloor$,
a formula which in fact matches all the best arrangements known so far except $n=5$ and $9$.
It would be nice to know if this construction really does solve the problem!


\section{Circles in the Plane}\label{SecWild}

\begin{figure}[!ht]
\centerline{\includegraphics[angle=0, width=4.0in]{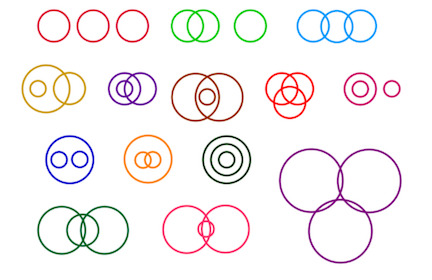}}
\caption{The fourteen ways to draw three circles in the affine plane.} \label{FigWildJG}
\end{figure}

\begin{figure}[!ht]
\centerline{\includegraphics[angle=0, width=6.0in]{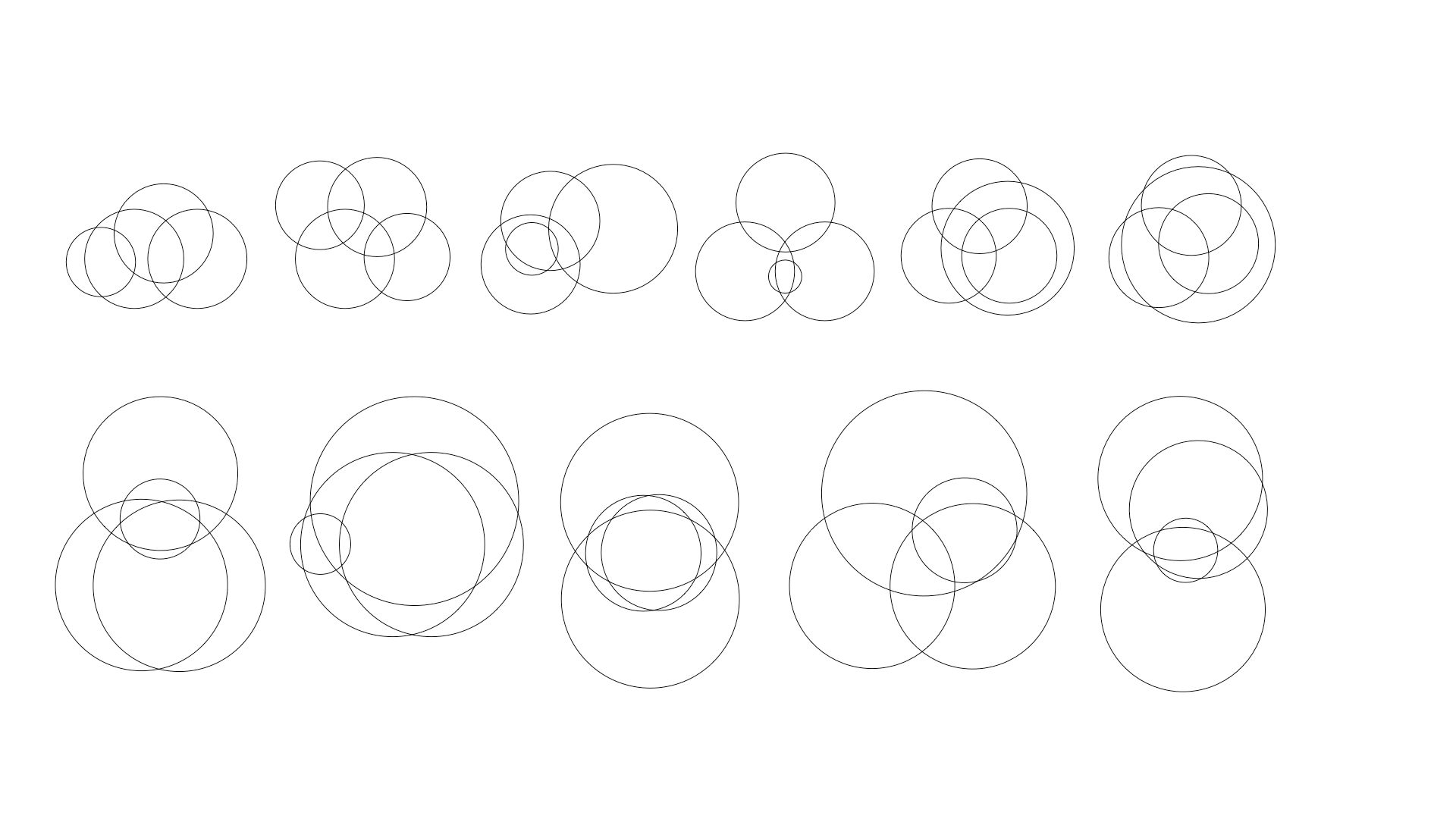}}
\caption{Eight of the $173$ ways to draw four circles. For the full set of $173$ drawings, see \seqnum{A250001}.}
\label{FigWild4a}
\end{figure}

\begin{figure}[!ht]
\centerline{\includegraphics[angle=0, width=6.0in]{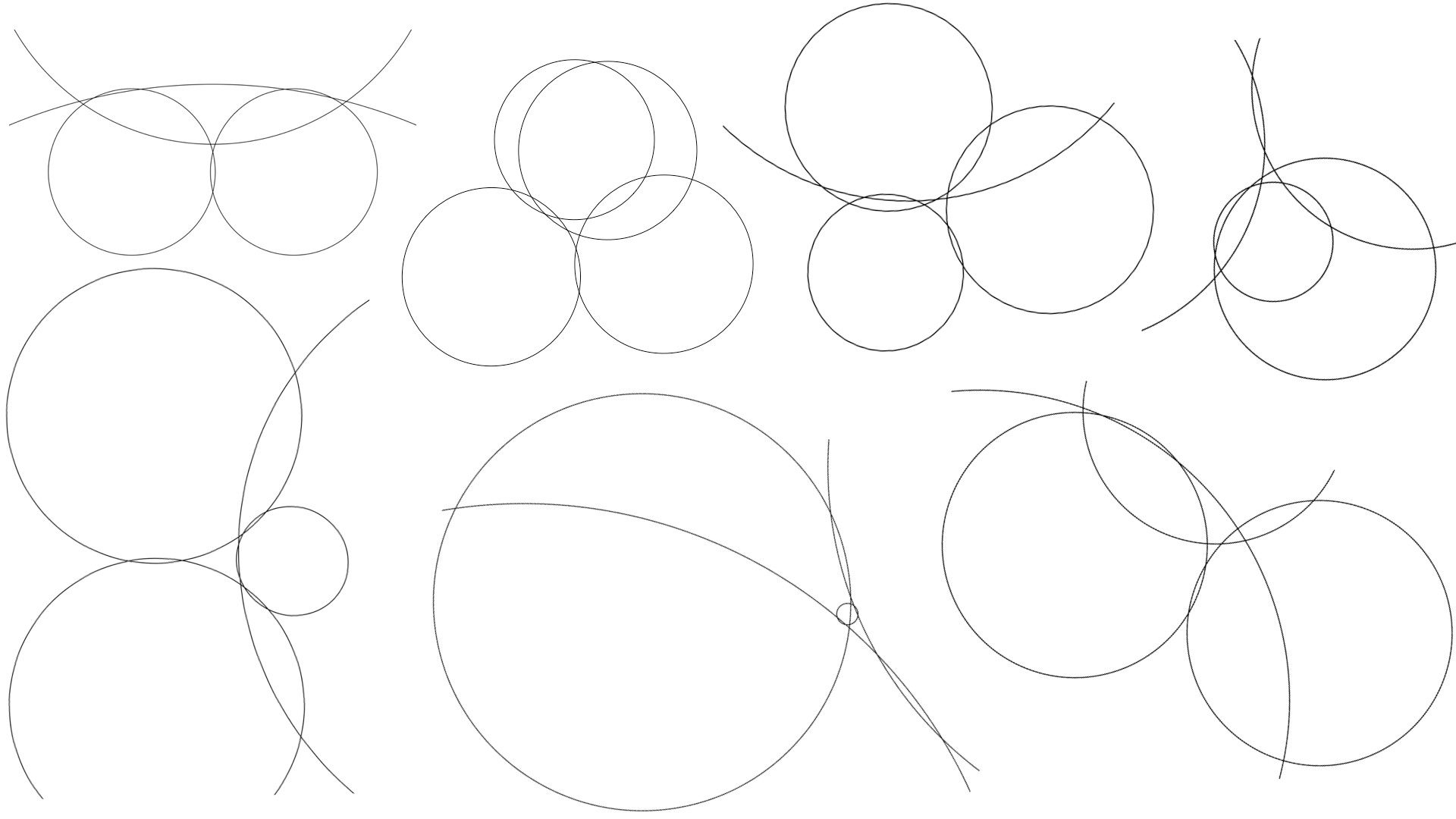}}
\caption{Seven further ways (out of 173) to draw four circles.}
\label{FigWild4b}
\end{figure}

The runner-up in the competition for A250000 is now
\seqnum{A250001}: here $a(n)$ is the number of ways to draw $n$ circles
in the affine plane. Two circles must be disjoint or meet in two distinct points (tangential 
contacts are not permitted), and
three circles may not meet at a point.\footnote{The circles may have different radii.  Two arrangements are considered the same if one can be continuously changed to the other while keeping all circles circular (although the radii may be continuously changed), without changing the multiplicity of intersection points, and without a circle passing through an intersection point. Turning the whole configuration over is allowed.}
The sequence was proposed by Jonathan Wild, a professor of music at
McGill University, who found the values 
$a(1)=1$, $a(2)=3$, $a(3)=14$, $a(4)=173$, and, jointly with Christopher  Jones, $a(5)=16951$
(see Figures~\ref{FigWildJG}-\ref{FigWild4b}).

\begin{figure}[!ht]
\centerline{\includegraphics[angle=270, width=3.0in]{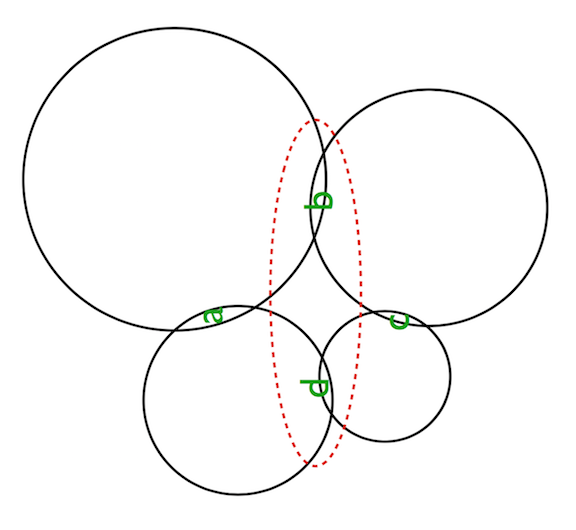}}
\caption{A hypothetical arrangement of five circles that can only be realized if
one or more of the circles is distorted. }
\label{FigWild5a}
\end{figure}

Wild and  Jones have
 found that there are complications which first appear
when five circles are being considered: here there are arrangements which theoretically
could exist if one considered only the intersections
between circles, but which cannot actually be drawn using circles.
For example, start with four circles arranged in a chain, each one overlapping its two neighbors, and label the overlaps a, b, c, d (see Figure~\ref{FigWild5a}). 
Suppose we try to add a fifth circle that meets all four circles but avoids their overlaps, encloses overlaps b and d, but does not enclose overlaps a or c. 
This can be drawn if the fifth circle is flattened to an ellipse, but
it can be shown that the arrangement cannot be realized with five circles.
There are $26$ such unrealizable arrangements of five circles,
which can be ruled out by ad hoc arguments. 

The delicate configurations like those in Figure~\ref{FigWild4b}
are very appealing.
It would be interesting to see all $17142$ arrangements of 
five or fewer circles displayed along the Great Wall of China.


\section{Earliest Cube-Free Binary Sequence}\label{SecCFBS}

There is an obvious way to sort integer sequences $a(1), a(2), a(3), a(4), ...$ into lexicographic order.
A number of recent entries in the OEIS are defined to be the
lexicographically earliest sequence of nonnegative or positive integers satisfying certain
conditions.

\vspace*{+.1in}

For example, one of the first results in the subject
now called   ``Combinatorics on Words''
was Axel Thue's 1912 theorem that  the ``Thue-Morse
sequence''  
$$
T ~=~ 0, 1, 1, 0, 1, 0, 0, 1, 1, 0, 0, 1, 0, 1, 1, 0, 1, 0, 0, 1, 0, 1, 1, 0, 0, 1,  \ldots
$$
(\seqnum{A010060}) contains no substring of the form $XXX$,
that is, $T$ is {\em cube-free}. $T$ can be defined as a fixed point
of the mapping $0  \to 01, 1 \to 10$; alternatively,
by taking $a(n)$ to be the parity of the number of $1$s in the binary expansion of $n$.
$105$ years later, David W. Wilson asked for the
lexicographically earliest cube-free sequence of $0$s and $1$s.
Using a back-tracking algorithm, he found what appear to be the first
$10000$ terms, which begin
\beql{EqCFBS1}
0, 0, 1, 0, 0, 1, 0, 1, 0, 0, 1, 0, 0, 1, 1, 0, 0, 1, 0, 0, 1, 0, 1, 0, 0, 1,  \ldots.
\eeq
This is now \seqnum{A282317}.  

There is no difficulty in showing that the sequence exists.\footnote{Thanks to Jean-Paul Allouche for this argument.}
To see this, make the set $S$ of all infinite binary sequences 
$a = (a(1), a(2), \ldots)$ into a metric space by defining $d(a,b)$ to be $0$ if 
$a=b$, or $2^{-i}$ if $a$ and $b$ first differ
at position $i$.  This identifies $S$ with the Cantor set in $[0,1)$.
The subset $F \subset S$ of infinite cube-free sequences is nonempty and has an infimum $c$ 
say. It is easy to show that the complement $S \setminus F$, sequences that contain a cube,
is an open set in this topology, so $F$ is closed and $c\in F$.

So far only the first $999$ terms of \seqnum{A282317} have been verified to be correct
(by showing that there is at least one infinite cube-free sequence with that beginning). 
The rest of the $10000$ terms are only conjectural.
It would be nice to know more.
In particular, does this sequence have an alternative construction?
There is no apparent formula or recurrence, which seems surprising.

\begin{figure}[!ht]
\centerline{\includegraphics[angle=0, width=5in]{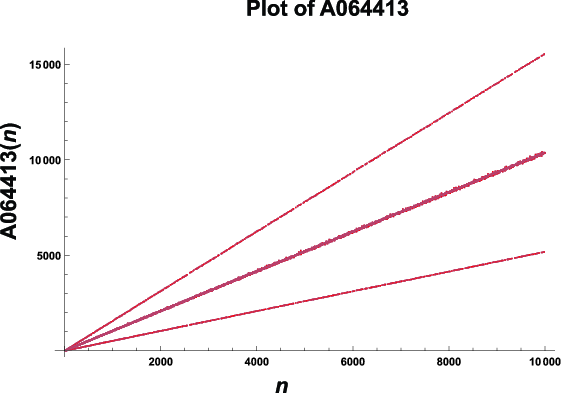}}
\caption{The first 10000 terms of the EKG sequence, so named because locally
this graph resembles an EKG.  Every number appears exactly once. }
\label{FigEKG1}
\end{figure}


\section{The EKG and Yellowstone Sequences}\label{SecEKG}

To continue the ``lexicographically earliest'' theme, many recent entries in the OEIS 
are defined to be the lexicographically earliest sequence $a(1), a(2), \ldots$ of 
distinct positive integers satisfying certain divisibility conditions.

The first task here is usually to show
that there are no missing numbers, i.e., that the sequence is a permutation
of the positive integers.
Sequences of this type were studied in a $1983$ paper by
Erd\H{o}s, Freud, and Hegyv\'{a}ri, 
which included the examples 
\seqnum{A036552} ($a(2n) =$ smallest missing number, $a(2n+1) = 2a(2n)$) and
\seqnum{A064736} ($a(2n+2) =$ smallest missing number, $a(2n+1) = a(2n) \cdot a(2n+2)$).
For these two it is clear that there are no missing numbers.  This is less obvious, but still true, for
Jonathan Ayres's {\em EKG sequence}, \seqnum{A064413}, defined to
be the lexicographically earliest sequence of distinct positive integers such that
$$
\gcd(a(n-1),a(n)) ~>~ 1 \mbox{~for~all~} n \ge 3\,.
$$
This begins
$$
1, 2, 4, 6, 3, 9, 12, 8, 10, 5, 15, 18, 14, 7, 21, 24, 16, 20, 22, 11, 33, 27, \ldots\,.
$$
The proof that it is a permutation is omitted--it is similar to the proof for the
Yellowstone sequence given below.

Next, one can investigate  the rate of growth.  
In the case of \seqnum{A064413}, the points appear to lie roughly on
three curved lines (Figure~\ref{FigEKG1}), although the following
conjecture of Lagarias, Rains, and Sloane (2002) is still open.

\begin{conj}\label{ConjEKG} In the 
EKG sequence \seqnum{A064413}, if $a(n)$ is neither a prime nor
three times a prime then
$$
a(n)  ~\sim~ n \left( 1 ~+~ \frac{1}{3 \log n}\right)\,;
$$
if $a(n)$ is  a prime  then
$$
a(n)  ~\sim~ \frac{1}{2} \,n \left( 1 ~+~ \frac{1}{3 \log n}\right)\,;
$$
and if $a(n)$ is $3$ times a prime  then
$$
a(n)  ~\sim~ \frac{3}{2} \, n \left( 1 ~+~ \frac{1}{3 \log n}\right)\,.
$$
\end{conj}

Furthermore, if the sequence is a permutation, one can also try to study its cycle
structure.  However, this often leads to very difficult
questions, similar to those encountered in studying the Collatz
conjecture,
and we can't do much more than collect experimental data.
Typically there is a set of finite cycles, and one or more apparently infinite cycles,
but we can't prove that the apparently
infinite cycles really are infinite, nor that they are distinct. See the entries for
\seqnum{A064413} and \seqnum{A098550} for examples.

\begin{figure}[!ht]
\centerline{\includegraphics[angle=0, width=4in]{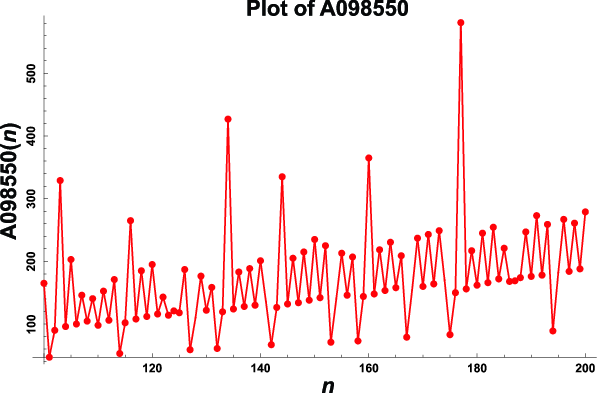}}
\caption{Plot of terms $a(101)$ through $a(200)$ of the Yellowstone sequence.
The sequence has a downward spike to $a(n)$ 
when $a(n)$ is a prime, and larger upward spikes
(the ``geysers'', which suggests the
name for this sequence) two steps later.}
\label{FigY1}
\end{figure}

\begin{figure}[!ht]
\centerline{\includegraphics[angle=0, width=6.75in]{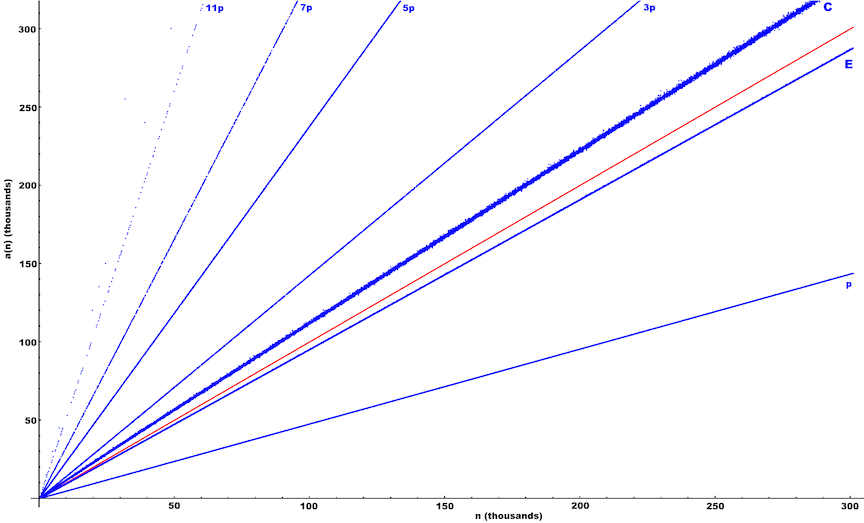}}
\caption{Scatterplot of the first 300,000 terms of the Yellowstone sequence. The primes lie
on the lowest line (labeled ``p''), the even numbers on the second line (``E''), 
the majority of the odd
composite numbers on the third line (``C''), and the 
$3p$, $5p$, $7p$, $11p, \ldots$
points on the higher lines. The lines are not actually straight,
except for the red line $f(x)=x$, which
is included for reference.}
\label{FigY2}
\end{figure}

The definition of the {\em Yellowstone sequence} 
(Reinhard Zumkeller, $2004$, \seqnum{A098550},  \cite{AHS15})
is similar to that of the EKG sequence, but now the requirement is that, for $n>3$,
$$
\gcd(a(n-2),a(n)) ~>~ 1 \mbox{~and~}  \gcd(a(n-1),a(n)) = 1\,.
$$
This begins
$$
1, 2, 3, 4, 9, 8, 15, 14, 5, 6, 25, 12, 35, 16, 7, 10, 21, 20, 27, 22, 39, 11, \ldots\,.
$$
Figure~\ref{FigY1} shows terms $a(101)=47$ through $a(200)=279$, with
successive points joined by lines.

\begin{thm}\label{ThY}
The Yellowstone sequence  \seqnum{A098550} is a permutation
of the positive integers.
\end{thm}

The proof is typical of the arguments used to
prove that several similar sequences
are permutations, including the EKG sequence above.
\begin{proof}
There are several steps.

\noindent (i) The sequence is infinite. (For $p\,a(n-2)$ is always a candidate for $a(n)$,
where $p$ is a prime larger than any divisor of $a(i)$, $i<n$.)

\noindent (ii) There are infinitely many different primes
that divide the terms of the sequence.
(If not, there is a prime $p$ such that all terms
are products of primes less than $p$. Using (i),
find a term $a(n) > p^2$, and let $q$ be a common prime factor of $a(n-2)$ and $a(n)$.
But now $pq<p^2 < a(n)$ is a smaller candidate for $a(n)$, a contradiction.)

\noindent (iii) For any prime $p$, some term is divisible by $p$. (For if not, no
prime $q>p$  can divide any $a(n)$: if $a(n)=kq$ is the first multiple of $q$
to appear, $kp$ would be a smaller candidate for $a(n)$. This  contradicts (ii).)

\noindent (iv) For any prime $p$, $p$ divides infinitely
many terms.  (If not, let $p^i$ be larger than any multiple of $p$
in the sequence, and choose a prime $q>p^i$. Again we obtain a contradiction.)

\noindent (v) Every prime $p$ is a term in the sequence. (Suppose not,
and using (i), choose $n_0$ such that $a(n)>p$ for all
$n > n_0$. Using (iv), find $a(n)=kp, k>1$, for
some $n>n_0$. But then $a(n+2)=p$, a contradiction.)

\noindent (vi) All numbers appear. For if not, let $k$ be the smallest missing
number, and choose $n_0$ so that all of $1, \ldots, k-1$ have occurred in
$a(1), \dots, a(n_0)$.
Let $p$ be a prime dividing $k$. Since, by (iv),
$p$ divides infinitely many terms,
there is a number $n_1 > n_0$ such that $\gcd(a(n_1),k) > 1$.
This forces
\beql{Eq1}
\gcd(a(n),k) > 1 \mbox{~for~\textbf{all}~}  n \ge n_1.
\eeq
(If not, there would be some $j \ge n_1$ where
$\gcd(a(j),k)>1$ and $\gcd(a(j+1),k)=1$, which would
lead to $a(j+2)=k$.) But \eqn{Eq1} is impossible, because we know from (v)
that infinitely many of the $a(n)$ are primes.
\end{proof}

The growth of this sequence is more
complicated than that of the EKG sequence.
Figure~\ref{FigY2} shows the first 300,000 terms, without lines connecting the points.
The points appear to fall on or close to a number
of distinct curves.  There is a conjecture in \cite[p.~5]{AHS15} that would explain these curves.


\section{Three Further Lexicographically Earliest Sequences}\label{SecLES3}
Here are three further examples of this type, all of which are 
surely permutations of the positive integers. 
For the first there is a proof, for the second there is ``almost'' a proof,
but the third may be beyond reach.

The first (Leroy Quet, $2007$, \seqnum{A127202}) is
the lexicographically earliest sequence of distinct positive integers 
such that
$$
\gcd(a(n-1), a(n)) ~\ne~ \gcd(a(n-2), a(n-1)) \quad  \mbox{for~} n \ge 3\,.
$$
It begins
$$
1, 2, 4, 3, 6, 5, 10, 7, 14, 8, 9, 12, 11, 22, 13, 26, 15, 18, 16, 17, 34, 19,  \ldots\,.
$$
For the second (R\'{e}my Sigrist,  2017, \seqnum{A280864}),
the definition is
$$
\mbox{if~a~prime~}p \mbox{~divides~} a(n),
\mbox{~then~it~divides~exactly~one~of~}
a(n-1) \mbox{~and~} a(n+1), \mbox{~for~} n \ge 2\,,
$$
and the initial terms are
$$
1, 2, 4, 3, 6, 8, 5, 10, 12, 9, 7, 14, 16, 11, 22, 18, 15, 20, 24, 21, 28, 26,  \ldots\,.
$$
The proof that the first is a permutation is similar to that for
the Yellowstone sequence, although a bit more involved
(see \seqnum{A127202}).
The second struck me as one of those ``drop everything and
work on this'' problems that are common hazards when editing new submissions to the OEIS.
However, after several months,  I could prove that every prime and every even number appears,
and that if $p$ is an odd prime then there are infinitely many odd multiples of $p$ 
(see \seqnum{A280864} for details), but I could not prove that every odd number appears.  The missing 
step feels like it is only a couple of cups of coffee away, and I'm hoping that 
some reader of this article will complete the proof.

The third example (Henry Bottomley,
2000, \seqnum{A055265}) is the lexicographically earliest sequence
of distinct positive integers such that $a(n-1)+a(n)$ is a prime for $n\ge 2$:
$$
1, 2, 3, 4, 7, 6, 5, 8, 9, 10, 13, 16, 15, 14, 17, 12, 11, 18, 19, 22, 21, 20, \ldots,.
$$
The terms appear to lie on or near the line $a(n)=n$, but the proof that every number appears
may be difficult because it involves the gaps between the primes.


\section{Two-Dimensional Lexicographically Earliest Arrays}\label{Sec2DLES}

\begin{figure}[!ht]
\centerline{\includegraphics[angle=0, width=3in]{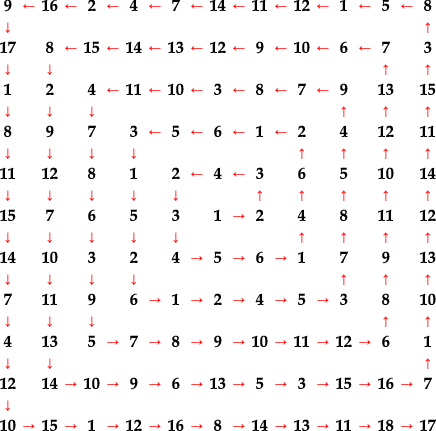}}
\caption{\seqnum{A274640}: choose smallest positive number so
that no row, column, or diagonal contains a repeat. Are the
rows, columns, diagonals permutations of $\NN$?
}
\label{FigSpiral}
\end{figure}

The OEIS is primarily a database of sequences $(a_n, n \ge n_0)$. 
However, triangles of numbers are included by reading them by rows.
Pascal's triangle becomes $1, ~1,1,~ 1,2,1,~ 1,3,3,1,$ $1,4,6,4,1,~ \ldots$,
which (without the extra spaces) is \seqnum{A007318}. Doubly-indexed arrays
$(T_{m,n}$,  $m \ge m_0, n \ge n_0)$ are converted to sequences by reading them
by antidiagonals (in either the upwards or downwards directions, or both). So an array $(T_{m,n},  m \ge 0, n \ge 0)$  
might become $T_{0,0}, T_{1,0}, T_{0,1},
T_{2,0}, T_{1.1}, T_{0,2}, \ldots$. For example, the table of 
Nim-sums $m \oplus n$:
\beql{EqNimSum}
\begin{array}{lllllllll}
0 & 1 & 2 & 3 & 4 & 5 & 6 & 7 & \cdots \\
1 & 0 & 3 & 2 & 5 & 4 & 7 & 6 & \cdots \\
2 & 3 & 0 & 1 & 6 & 7 & 4 & 5 & \cdots \\
3 & 2 & 1 & 0 & 7 & 6 & 5 & 4 & \cdots \\
4 & 5 & 6 & 7 & 0 & 1 & 2 & 3 & \cdots \\
5 & 4 & 7 & 6 & 1 & 0 & 3 & 2 & \cdots \\
6 & 7 & 4 & 5 & 2 & 3 & 0 & 1 & \cdots \\
7 & 6 & 5 & 4 & 3 & 2 & 1 & 0 & \cdots \\
\cdot & \cdot & \cdot & \cdot & \cdot & \cdot & \cdot & \cdot & \cdots
\end{array}
\eeq
produces the sequence \seqnum{A003987}:
$$
0,~ 1, 1,~ 2, 0, 2,~ 3, 3, 3, 3,~ 4, 2, 0, 2, 4,~ 5, 5, 1, 1, 5, 5,~ 6, 4, 6, 0, 6, 4, 6,~  \ldots\,. 
$$
Doubly-indexed doubly-infinite arrays $(T_{m,n}, m \in \ZZ, n \in \ZZ)$ can become sequences
by reading them in a spiral around the origin, in say a counter-clockwise
direction:
$T_{0,0}, T_{1,0}, T_{1,1}, T_{0,1},$ $T_{-1,1},
T_{-1,0}, T_{-1,-1}, T_{0,-1}, \ldots$
(cf.~Figure~\ref{FigSpiral}).

There are many ``lexicographically earliest'' versions of these arrays.
For example, the Nim-sum array \eqn{EqNimSum} has an equivalent definition: 
scan  along upwards antidiagonals, filling in each cell with the
smallest nonnegative number that is neither in the row
to the left of that cell nor in the column above it. 

A variation on the Nim-sum array was proposed by Alec Jones in 2016, as a kind of 
``infinite Sudoku array''. This array $(T_{m,n},  m \ge 0, n \ge 0)$  is to be filled in
by upwards antidiagonals, always choosing the smallest positive integer such that
no row, column, diagonal, or antidiagonal contains a repeated term.
The top left corner of the array is:

\beql{EqAJ1} 
\begin{array}{ccccccccc}
1 &  3 &  2 &  6 &  4 &  5 &  10 &  11 &  \cdots \\
2 &  4 &  5 &  1 &  8 &  3 &   6 &  12 &  \cdots \\
3 &  1 &  6 &  2 &  9 &  7 &   5 &  4 &  \cdots \\
4 &  2 &  3 &  5 &  1 &  8 &   9 &  7 &  \cdots \\
5 &  7 &  1 &  4 &  2 &  6 &   3 &  15 &  \cdots \\
6 &  8 &  9 &  7 &  5 &  10 &  4 &  16 &  \cdots \\
7 &  5 &  4 &  3 &  6 &  14 &  8 &  9 &  \cdots \\
8 &  6 &  7 &  9 &  11 &  4 &  13 &  3 &  \cdots \\
\cdot & \cdot & \cdot & \cdot & \cdot & \cdot & \cdot & \cdot & \cdots
\end{array} 
\eeq

\noindent
The resulting sequence (\seqnum{A269526}) is
$$
1,~ 2, 3,~ 3, 4, 2,~ 4, 1, 5, 6,~ 5, 2, 6, 1, 4,~ 6, 7, 3, 2, 8, 5,~ 7, 8, 1, 5, 9, 3, 10,~  \ldots\,.
$$

This array has many interesting properties.
If we subtract $1$ from each entry, the entries are the Nim-values for a game
played with two piles of counters, of sizes $m$ and $n$, and reminiscent
of Wythoff's game (see \seqnum{A004481}, \seqnum{A274528}).

But the main question about the array \eqn{EqAJ1} is,
are the individual rows, columns, and diagonals of this array 
permutations of $\NN$?
(The antidiagonals are obviously not, since they are finite sequences.)
It is easy to see that each column is a permutation. In column $c \ge 0$, a number $k$ 
will eventually be the smallest missing number and will
appear in some cell in that column, unless there is a copy of $k$ to the 
North-West, West, or South-West of that cell.
But there are at most $c$ copies of $k$ in all the earlier columns, 
so eventually $k$ {\em will} appear. 

The rows are also permutations, although the proof is less
obvious. Consider row $r \ge 0$, and suppose $k$ never appears. 
There are at most $r$ copies of $k$ in the earlier rows, and
these can affect only a bounded portion of row $r$.
Consider a cell $(r,n)$, $n \ge 0$ large. If $k$ is not to appear in that cell, there must be 
a copy of $k$ in the antidiagonal to the South-West.
So in the triangle bounded by row $r$, column $0$, and the 
antidiagonal through $(r,n)$, there must be at least $n+1-r$ copies of $k$.
Imagine these $k$s replaced by chess queens. By construction they are
mutually non-attacking. But it is known (\cite[Problem~$252$]{VGL02},
or \seqnum{A274616})
that on a triangular  half-chessboard of side $n$,
there can be at most $2n/3+1$ mutually non-attacking queens, which
for large $n$ leads to a contradiction. 

As to the diagonals, although they appear to be permutations, 
this is an open question. The argument using non-attacking queens breaks down because 
the diagonal of the half-chessboard contains only half as many squares as the sides.
Even the main diagonal, \seqnum{A274318},
$$
1, 4, 6, 5, 2, 10, 8, 3, 7, 9, 16, 26, 29, 22, 20, 23, 28, 38, 12, 32, 46, 13, 14, 11, 15, \ldots\,,
$$
is not presently known to be a permutation of $\NN$.

\begin{figure}[!ht]
\centerline{\includegraphics[angle=0, width=5in]{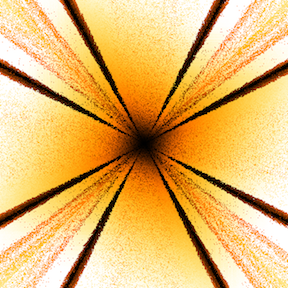}}
\caption{Colored representation of central $200 \times 200$ portion of the spiral in Figure~\ref{FigSpiral}: the colors represent the values, ranging from black (smallest) to white (largest).}
\label{FigSpiralColor}
\end{figure}

The spiral version of this array is even
more frustrating. This array ($(T(m,n), m \in \ZZ, n \in \ZZ)$, \seqnum{A274640},
proposed by Zak Seidov and Kerry Mitchell in June $2016$), is constructed in a counterclockwise spiral,
filling in each cell with the smallest positive number
such that no row, column, or diagonal 
contains a repeated term (Figures~\ref{FigSpiral}, \ref{FigSpiralColor}).
(``Diagonal'' now means any line of cells of slope $\pm 1$.)

Although it seems very plausible that every row, column, and diagonal 
is a permutation of $\NN$, now there are no proofs at all. The eight 
spokes through the center are sequences \seqnum{A274924}--\seqnum{A274931}.
For example, the row through the central cell is
$$
\ldots, 14, 25, 13, 17, 10, 15, 7, 6, 5, 3, 1, 2, 4,  8, 11, 12, 16,  9, 19, 24, 22, \ldots\,,
$$
which is \seqnum{A274928} reversed followed by \seqnum{A274924}.
Is it a permutation of $\NN$? We do not know.


\section{Fun With Digits}\label{SecDavis}

Functions of the digits of numbers have always 
fascinated people,\footnote{Although in {\em A Mathematician's Apology}, G. H. Hardy, 
referring to the fact that 
$1089$ and $2178$ are the smallest numbers which when
written backwards are nontrivial multiples of themselves (cf. \seqnum{A008919}), remarked
that this fact was ``likely to amuse amateurs'', but was not of interest to mathematicians.}
and one such function was in the news in 2017. The idea underlying this story
and several related sequences is to start with some simple function $f(n)$
of the digits of $n$ in some base, iterate it, and watch what happens.

For the first example we write $n$ as a product of prime powers,
$n = p_1^{e_1} p_2^{e_2} \cdots$ with the $p_i$ in increasing order, and define
$f(n)$ to be the decimal concatenation $p_1 e_1 p_2 e_2 \ldots$,
where we omit any exponents $e_i$ that are equal to $1$. 
So $f(7)=f(7^1)=7, f(8)=f(2^3)=23$.

The initial values of $f(n)$ (\seqnum{A080670}) are
$$
\begin{array}{crrrrrrrrrrrrrrrrrrrrr}
n: & 1 & 2 & 3 & 4 & 5 & 6 & 7 & 8 & 9 & 10 & 11 & 12 & 13 & 14 & 15 & 16 & 17 & 18 & \ldots \\
f(n): & 1 & 2 & 3 & 22 & 5 & 23 & 7 & 23 & 32 & 25 & 11 & 223 & 13 & 27 & 35 & 24 & 17 & 232 & \ldots
\end{array}
$$
If we start with a positive number $n$ and repeatedly apply $f$, in many small cases
we rapidly reach a prime (or $1$).\footnote{In what follows
we will tacitly assume $n \ge 2$, to avoid having
to repeatedly say ``(or $1$)''.}
For example, $9 = 3^2 \to 32 = 2^5 \to 25 = 5^2 \to 52 = 2^2 13 \to 2213$,
a prime.
Define $F(n)$ to be the prime that is eventually reached,
or $-1$ if the iteration never reaches a prime. 
The value $-1$ will occur if the iterates are unbounded or if
they enter a cycle of composite numbers.
The initial values of $F(n)$ (\seqnum{A195264}) are
$$
1, 2, 3, 211, 5, 23, 7, 23, 2213, 2213, 11, 223, 13, 311, 1129, 233, 17, 17137, 19, \ldots\,.
$$
$F(20)$ is currently unknown (after $110$ steps the trajectory of $20$ 
has stalled at a $192$-digit number which has not yet been factored).
At a DIMACS conference in October $2014$ to celebrate the $50$th anniversary of the start
of what is now the OEIS, John H. Conway offered \$1000 for
a proof or disproof of his conjecture that the iteration of $f$ will always
reach a prime.

However, in June 2017 James Davis found a number $D_0 = 13532385396179$
whose prime factorization is $13 \cdot 53^2 \cdot 3853 \cdot 96179$,
and so clearly $f(D_0) = D_0$ and $F(D_0)=-1$.

The method used by James Davis to find $D_0$ is quite simple. Suppose $n = m \cdot p$ is fixed by $f$, where $p$ is a prime greater than
all the prime factors of $m$. Then $f(n)=f(m)10^y+p$, where $y$ is
the number of digits in $p$. From $f(n)=n$ we have
$p = \tfrac{f(m) 10^y}{m-1}$. 
Assuming $p \ne 2, 5$, this implies that $p$ divides $f(m)$, and setting $x=f(m)/p$,
we find that $m=x 10^y+1$ with  $p = \tfrac{f(m)}{x}$ prime.
A computer easily finds the solution
$x=1407$, $y=5$, $m=140700001$, $p=96179$, and so $n=D_0$.

No other composite fixed points are known, and David J. Seal has recently shown
that there is no composite fixed point less than $D_0$.
It is easy, however, to find numbers whose trajectory under $f$ 
ends at $D_0$, by repeatedly finding a prime prefix of the previous number,
as shown by the example\footnote{Found by Hans Havermann.}
$D_1 = 13^{532385396179}$ with $f(D_1) = D_0$.
So presumably there are infinitely many $n$ with $F(n)=-1$.

Consideration of the analogous questions in other bases might have
suggested that counterexamples to Conway's question could exist.
We will use subscripts to indicate the base (so $4_{10} = 100_2$).
The base-$2$ analog of $f$, $f_2$ (say), is defined by taking
$f_2(p_1^{e_1} p_2^{e_2} \cdots)$ to be the concatenation $p_1 e_1 p_2 e_2 \ldots$, as before
(again omitting any $e_i$ that are $1$),
except that now we write the $p_i$ and $e_i$ in base $2$ and interpret the concatenation
as a base-$2$ number. For example, $f_2(8) = f_2(2^3) = 1011_2 = 11_{10}$.

The initial values of $f_2(n)$ (\seqnum{A230625}) are
$$
\begin{array}{crrrrrrrrrrrrrrrrrrrrr}
n: & 1 & 2 & 3 & 4 & 5 & 6 & 7 & 8 & 9 & 10 & 11 & 12 & 13 & 14 & 15 & 16 & 17 & 18 & \ldots \\
f_2(n): & 1 & 2 & 3 & 10 & 5 & 11 & 7 & 11 & 14 & 21 & 11 & 43 & 13 & 23 & 29 & 20 & 17 & 46 &  \ldots
\end{array}
$$
and the base-$2$ analog of $F$, $F_2$ (\seqnum{A230627}) is the prime (or $1$)
that is reached when $f_2$ is repeatedly applied to $n$, or $-1$ if no prime (or $1$) is 
reached:
$$
1, 2, 3, 31, 5, 11, 7, 11, 23, 31, 11, 43, 13, 23, 29, 251, 17, 23,  \ldots\,.
$$
Now there is a fairly small composite fixed point, namely $255987$, found
by David J. Seal. Sean A. Irvine and Chai Wah Wu have also
studied this sequence, and the present status is that $F_2(n)$ is known
for all $n$ less than $12388$. All numbers in this range reach $1$, a prime,
the composite number $255987$, or one of the two cycles $1007 \longleftrightarrow 1269$ or
$1503 \longleftrightarrow 3751$. The numbers for which
$F_2(n)=-1$ are $217, 255, 446, 558, \ldots$ (\seqnum{A288847}).
Initially it appeared that $234$ might be
on this list, but Irvine found that after $104$ steps the trajectory reaches the 
$51$-digit prime 
$$
350743229748317519260857777660944018966290406786641\,.
$$


\section{Home Primes}\label{SecHome}
A rather older problem arises if we change the definition of $f(n)$ slightly, making
$f(8) = 222$ rather than $23$. 
So if $n = p_1 \cdot p_2 \cdot p_3 \cdot \ldots$, where
$p_1 \le p_2 \le p_3 \le \ldots$, then $f(n)$ is
the decimal concatenation $p_1 p_2 p_3 \ldots$ (\seqnum{A037276}).
In $1990$, Jeffrey Heleen studied the analog
of $F(n)$ for this function: that is, $F(n)$ is
the prime reached if we start with $n$ and repeatedly apply $f$, or $-1$
if no prime is ever reached (\seqnum{A037274}).

The trajectory of $8$ now takes $14$ steps to reach a prime (the individual prime factors here have been separated by spaces):
\begin{align} 
8 &  \to
2 ~ 2 ~ 2 \to
2 ~ 3 ~ 37 \to
3 ~ 19 ~ 41 \to
3 ~ 3 ~ 3 ~ 7 ~ 13 ~ 13 \to
3 ~ 11123771 \to
 7 ~ 149 ~ 317 ~ 941 \to \nonumber \\
& \to 229 ~ 31219729 \to
11 ~ 2084656339 \to
3 ~ 347 ~ 911 ~ 118189 \to
11 ~ 613 ~ 496501723 \to \nonumber \\
& \to 97 ~ 130517 ~ 917327 \to
53 ~ 1832651281459 \to
3 ~ 3 ~ 3 ~ 11 ~ 139 ~ 653 ~ 3863 ~ 5107 \nonumber \\
& \to 3331113965338635107\,, \nonumber
\end{align}
the last number being a prime.

Since $f(n)>n$ if $n$ is composite, now there cannot be any composite
fixed points nor any cycles of length greater than $1$. The only
way for $F(n)$ to be $-1$ is for the trajectory
of $n$ to be unbounded. This appears to be a harder problem
than the one in the previous section,
since so far no trajectory has been proved to be unbounded.
The first open case is $n=49$, which after $119$ iterations has reached a $251$-digit composite number (see \seqnum{A056938}).
The completion of the factorization for step $117$ took $765$ days by the general number field sieve, and at the time (December 2014) was one of the hardest factorizations ever completed.


\section{Power Trains}\label{SecPowerTrains}

A third choice for $f(n)$ was proposed by John H. Conway in 
$2007$: he called it the {\em power train} map.
If the decimal expansion of $n$ is $d_1 d_2 d_3 \ldots d_k$
(with $0 \le d_i \le 9$, $0 < d_1$), then
$f(n) = d_1^{d_2} \cdot d_3^{d_4}  \cdots$, ending 
with $\ldots \cdot d_k$ if $k$ is odd, or with $\ldots \cdot d_{k-1}^{d_k}$ 
if $k$ is even (\seqnum{A133500})..
We take $0^0$ to be $1$.
For example,
$f(39) = 3^9 = 19683$,
$f(623) = 6^2 \cdot 3 = 108$.
Conway observed that $2592 = 2^5 9^2$ is a non-trivial fixed point, and asked me if there were any others.
I found one more: 
$n = 2^{46} \cdot 3^6 \cdot 5^{10} \cdot 7^2 = 24547284284866560000000000$,
for which $f(n) = 2^4 \cdot 5^4 \cdot 7^2 \cdot 8^4 \cdot 2^8 \cdot 4^8 \cdot 6^6 \cdot 5^6 \cdot 0^0 \cdot 0^0 \cdot 0^0 \cdot 0^0 \cdot 0^0 = n$.
The eleven known fixed points 
(including the trivial values $1, \ldots,9$)
form \seqnum{A135385}, and it is known
that there are no further terms below $10^{100}$.
Maybe this is a hint that  for all of the functions $f(n)$ that
have just been mentioned, there may be only a handful
of genuinely exceptional values?


\section{A Memorable Prime}\label{SecMemorable}
If you happen to need an explicit $20$-digit prime in a hurry, it is useful to remember that although
$1, 121=11^2, 12321=111^2, 1234321=1111^2,  \ldots,$ and $12345678987654321 = 111111111^2$ are
not primes, the next term in \seqnum{A173426} {\em is} a prime,
$$
12345678910987654321\,.
$$
As David Broadhurst remarked on the Number Theory Mailing List
in August 2015, this is a memorable prime!  He also pointed out that on
probabilistic grounds, there should be infinitely many values of $n$ such
that the decimal concatenation of the  numbers $1$ up through $n$ followed by $n-1$ down through $1$
is a prime. Shortly afterwards, Shyam Sunder Gupta found what is presumably
the next prime in the sequence, corresponding to $n=2446$, the
$17350$-digit probable prime $1234567..244524462445..7654321$.
Serge Batalov has shown that there are no further terms with $n<60000$.
What is the next term? The values $10, 2446$ are not enough to create an OEIS entry.


\section{A Missing Prime}\label{SecMissing}
The previous question naturally led me to wonder what the first prime is in the
simpler sequence (\seqnum{A007908}):
$$
1,12,123,1234,\ldots, 12345678910, 1234567891011, \ldots\,,
$$
formed by the decimal concatenation of the numbers $1$ through $n$.
In {\em Unsolved Problems in Number Theory}, Richard K. Guy reports that this question
was already asked by Charles Nicol and John Selfridge.
However, although the same probabilistic argument suggests that
there should be an infinite number of
primes of this type, not a single one is known.  I asked several friends to
help with the search, and as a result this task was taken up by the
folks who run the GIMP (or Great Internet Mersenne Prime) search,
and there is now a 
web page\footnote{\url{http://mersenneforum.org/showthread.php?t=20527}.}
that shows the current status of the search for the first prime.  As of August 2017 the search seems to have stalled,
the present status being  that all the potential values of $n$ through $344869$ failed
(obviously many values of $n$ can be ruled out by congruence conditions).
In this range the candidates have about two million digits. 
One estimate suggests that there is a probability of about $0.5$
that a prime will be found with $n < 10^6$, so it
would be good to resume this search.


\section{Post's Tag System}\label{SecTag}

In his recent book {\em Elements of 
Mathematics: From Euclid to G{\"o}del}\footnote{A superb
successor to Felix Klein's 1908  {\em Elementary Mathematics
from an Advanced Standpoint.}}, John Stillwell mentions that
Emil L. Post's tag system from the 1930s is still not understood.
Post asked the following question.  Take a finite string, or word,  $S$ of $0$s and $1$s,
and if it begins with $0$,
append $00$ to the end of $S$ and delete the first three symbols, or
if it begins with $1$,
append $1101$ to the end of $S$ and delete the first three symbols.
When this process is iterated,
eventually one of three things will happen: either $S$ will reach
the empty word ($S$ {\em dies}), $S$ will enter a loop ($S$ {\em cycles}), or
$S$ will keep growing for ever ($S$ {\em blows up}).
For example, $S=1000$ reaches the empty word $\epsilon$ at the $7$th step:
$$
1000 \to  01101 \to  0100 \to  000 \to  00 \to  0 \to  \epsilon \,,
$$
whereas $100100$ enters a cycle of length six (indicated by
parentheses) after $15$ steps:
   \begin{align}\label{EqTag1}
100100 & \to  1001101 \to  11011101 \to  111011101 \to  0111011101 \to  101110100 \to  1101001101 \nonumber \\ 
&  \to  10011011101 \to 110111011101 \to  1110111011101 \to  01110111011101  \nonumber \\
& \to  1011101110100 \to  11011101001101
  \to  111010011011101 \to  0100110111011101 \nonumber \\
&  \to  (011011101110100 \to  01110111010000 \to  1011101000000 \to  11010000001101 \nonumber \\
  &  \to  100000011011101 \to  0000110111011101) \,.
   \end{align}
Post was hoping to find an algorithm which, given $S$,
would determine which of these outcomes would occur.
He did not succeed.

Post called this process a `tag system.' It can be generalized 
by considering initial words over an alphabet of size $M$ (rather than $2$),
allowing any fixed set $\sA$  of $M$ tag words to be appended (rather than $00$ and $1101$),
and deleting some fixed number $P$ of initial symbols at each step (not necessarily $3$).  
In 1961, Marvin Minsky showed that such a generalized tag system
could simulate a Turing machine. By choosing an appropriate alphabet,
an appropriate set $\sA$ of tag words to be appended, and an appropriate
value of $P$ (in fact $P=2$ will do), any computable function
can be simulated.
So, because of the undecidability of the Halting Problem, for general tag systems it is impossible to
predict which initial words will blow up. 

But what about Post's original tag system? Could this simulate a Turing machine
(by encoding the problem in the initial word $S$)?
At first this seems very unlikely, but the Cook-Wolfram theorem
 that the one-dimensional cellular automaton defined
by Rule $110$ can simulate a Turing machine (by encoding the
problem in the starting state) suggests that it might be possible.
If it {\em is} possible, there must be some initial words that blow up (again because
of the Halting Problem).

\begin{figure}[!ht]
\centerline{\includegraphics[angle=0, width=5.4in]{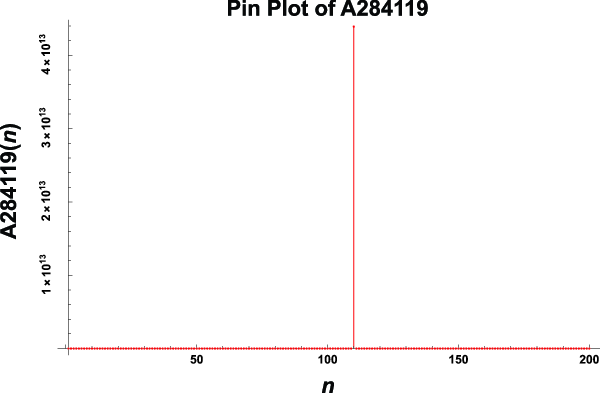}}
\caption{Pin plot illustrating Lars Blomberg's remarkable
discovery that the Post tag system started at the word $(100)^{110}$
takes an exceptionally long time ($43913328040672$  steps) to converge.}
\label{FigTag1}
\end{figure}

In early 2017, when I read Stillwell's book,, the OEIS contained
three sequences related to the original tag system, based on the work of
Peter Asveld  and submitted by Jeffrey Shallit:
\seqnum{A284116}, giving the maximal number of words in the
`trajectory' of any initial word $S$ of length $n$ ($18$ terms
were known), and two sequences connected with the especially
interesting starting word $\sigma_n$ of
length $3n$ consisting of $n$ copies of $100$.
\seqnum{A284119}$(n)$ is defined to be the number of words in the
trajectory of $\sigma_n$ before it enters a cycle or dies, or $-1$ if the trajectory blows up,
and \seqnum{A284121}$(n)$ is the length of the cycle, or $1$ if the trajectory dies, or $-1$
if the trajectory blows up.
For example, from \eqn{EqTag1} we see that \seqnum{A284119}$(2) = 15$ and
\seqnum{A284121}$(2) = 6$. Shallit had extended Asveld's work and had found $43$ terms 
of the two last-mentioned sequences.

I then added many further sequences based on tag systems 
discussed by Asveld, Liesbeth De~Mol,
Shigeru Watanabe, and others,
and appealed to contributors to the OEIS to extend them.

The most interesting response came from Lars Blomberg, who
investigated the trajectory of $\sigma_n$ for $ n \le 110$.
On September 9 2017 he reported that every $\sigma_n$ for $n \le 110$ had
either died or cycled after at most $13$ million terms,
except for $\sigma_{110}$, which after $38.10^{11}$ steps had reached a word of length 
$10^7$ and was still growing. This was exciting news!  Could $\sigma_{110}$ be
the first  word to be discovered that blew up?\footnote{Of course the fact that the same number $110$ was involved could not possibly be anything more than a coincidence.}
Sadly, on October 4 2017, Blomberg reported that after  $43913328040672$ 
steps $\sigma_{110}$ had terminated in the empty word.

\begin{figure}[!ht]
\centerline{\includegraphics[angle=0, width=6.0in]{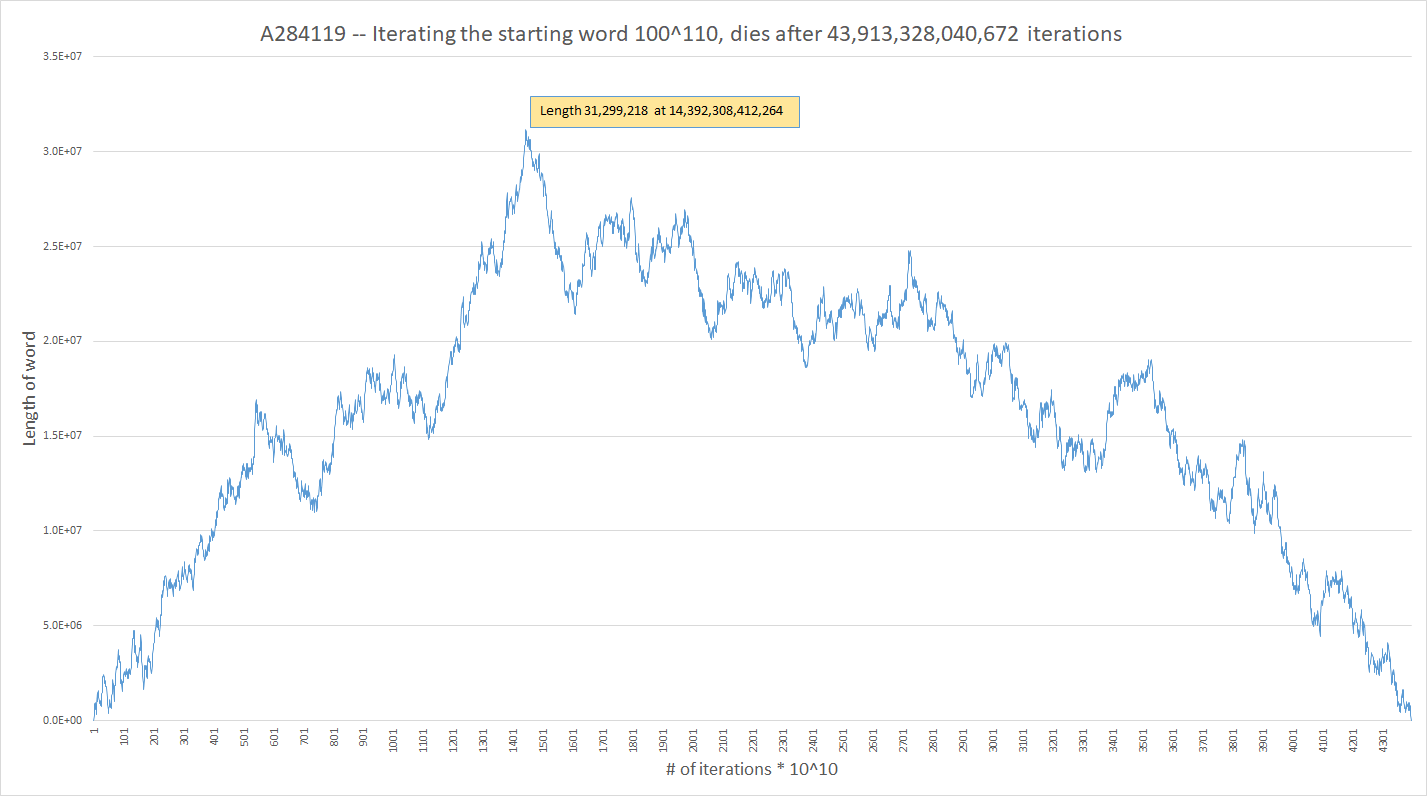}}
\caption{Lengths of successive words in trajectory of $(100)^{110}$
under the Post tag system. 
The numbers on the horizontal axis are spaced at multiples of $10^{12}$.}
\label{FigTag2}
\end{figure}

Figure~\ref{FigTag1} displays the remarkable graph (technically, a pin plot)
of the number of steps for $\sigma_n$ to either die or cycle for $n \le 200$.
Figure~\ref{FigTag2} shows the lengths of the
successive words in the trajectory of $\sigma_{110}$. 

In the past six months Blomberg has continued this investigation
and has determined the fate of $\sigma_n$ for all $n \le 6075$.
The new record-holder for the number of steps
before the trajectory dies is now held by $\sigma_{4974}$,
which takes $57042251906801$ steps, while $\sigma_{110}$
is in second place.

Of course it is still possible that some initial word $S$, not
necessarily of the form $\sigma_n$, will blow up, but this seems increasingly 
unlikely.  So Post's tag system probably does not simulate a Turing machine.

The question as to which $\sigma_n$ die and which cycle remains a mystery.
Up to $n=6075$, Blomberg's results show that about
one-sixth of the values of $n$ die and  five-sixths cycle. The precise
values can be found in \url{A291792}. It would be nice to understand this 
sequence better.


\section{Coordination Sequences}\label{SecCS}
This final section is concerned with {\em coordination sequences},
which arise in crystallography and in studying tiling problems, have
beautiful illustrations, and lead to many unsolved mathematical questions.

\begin{figure}[htb]
\centerline{\includegraphics[angle=0, width=.9\textwidth]{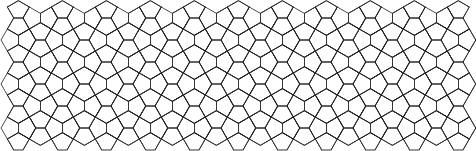}}
\caption{A portion of the Cairo tiling.} 
\label{FigCairo}
\end{figure}

The ``Cairo'' tiling, so called because it is said to be used
 on many streets in that city, is shown in Figure~\ref{FigCairo}.
Let $G$ denote the corresponding infinite graph (with vertices for points
where three or more tiles meet, and edges between two vertices where two tiles
meet). The figure is also a picture of the graph.

The distance 
between vertices $P, Q \in G$ is defined to be the number of
edges in the shortest path joining them. The 
{\em coordination sequence} of $G$ with
respect to a vertex $P \in G$ is then the sequence
$a(n)$ $(n \ge 0)$ giving the number of vertices $Q$ at distance $n$ from $P$.
Coordination sequences have been studied by crystallographers for 
many years \cite{OH1980}.

\begin{figure}[htb]
\centerline{\includegraphics[width=.35\textwidth]{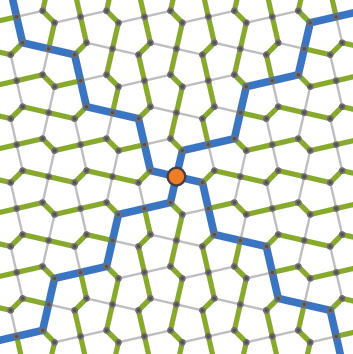}}
\caption{ The ``coloring book'' method applied to a tetravalent
vertex (the red dot) in the Cairo tiling, used to prove that
the coordination sequence is the same as that for the square grid.}
\label{FigTet2}
\end{figure}

The graph of the Cairo tiling has two kinds of vertices, trivalent (where three edges meet) and tetravalent.
As can be seen from  Figure~\ref{FigCairo}, the coordination sequence
with respect to a tetravalent vertex begins $1, 4, 8, 12, 16, 20, 24, \ldots$,
which appears to be the same as the coordination sequence \seqnum{A008574}
 for a vertex in the
familiar square grid. This observation seemed to be new.
Chaim Goodman-Strauss and I thought that such a simple 
fact should have a simple proof, and we developed an elementary
``coloring book''  procedure
\cite{GSS18} which not only proved this result but also established a number
of conjectured  formulas for coordination sequences of other tilings mentioned in entries in
the OEIS. 
The ``coloring book'' drawing of the Cairo graph of the Cairo graph centered
at a tetravalent vertex is shown in Figure~\ref{FigTet2}. This coloring
makes it easy to prove that
the coordination sequence is given by  $a(n)=4n$ for $n \ge 1$
(see  \cite{GSS18} for details).

For a trivalent vertex in the Cairo tiling, the coordination sequence is
$$
1, 3, 8, 12, 15, 20, 25, 28, 31, 36, 41, 44, 47, 52, 57, 60, 63, 68, \ldots
$$
(this is now \seqnum{A296368}), and we \cite{GSS18} show that
for $n \ge 3$, $a(n)=4n$ if $n$ is odd, $4n-1$ if $n \equiv 0 \pmod{4}$,
and $4n+1$ if $n \equiv 2 \pmod{4}$.

One can similarly define coordination sequences for other two- and higher-dimensional
structures, and the OEIS presently contains over $1300$ such sequences.
Many more could be added.
There are many excellent web sites with lists of tilings and crystals.
Brian Galebach's web site\footnote{\url{http://probabilitysports.com/tilings.html}.} 
is especially important, as it includes pictures of all ``$k$-uniform''
tilings with $k \le 6$, with over $1000$ tilings.
Darrah Chavey's article \cite{Chavey} and the  
Michael Hartley and Printable Paper web sites\footnote{
\url{http://www.dr-mikes-math-games-for-kids.com/archimedean-graph-paper.html}, \\
\url{https://www.printablepaper.net/category/graph}.}
have many further pictures, and the RCSR  and ToposPro
databases\footnote{\url{http://rcsr.net},  \url{http://topospro.com}.}
have thousands more.

\begin{figure}[htb]
\centerline{\includegraphics[width=.75\textwidth]{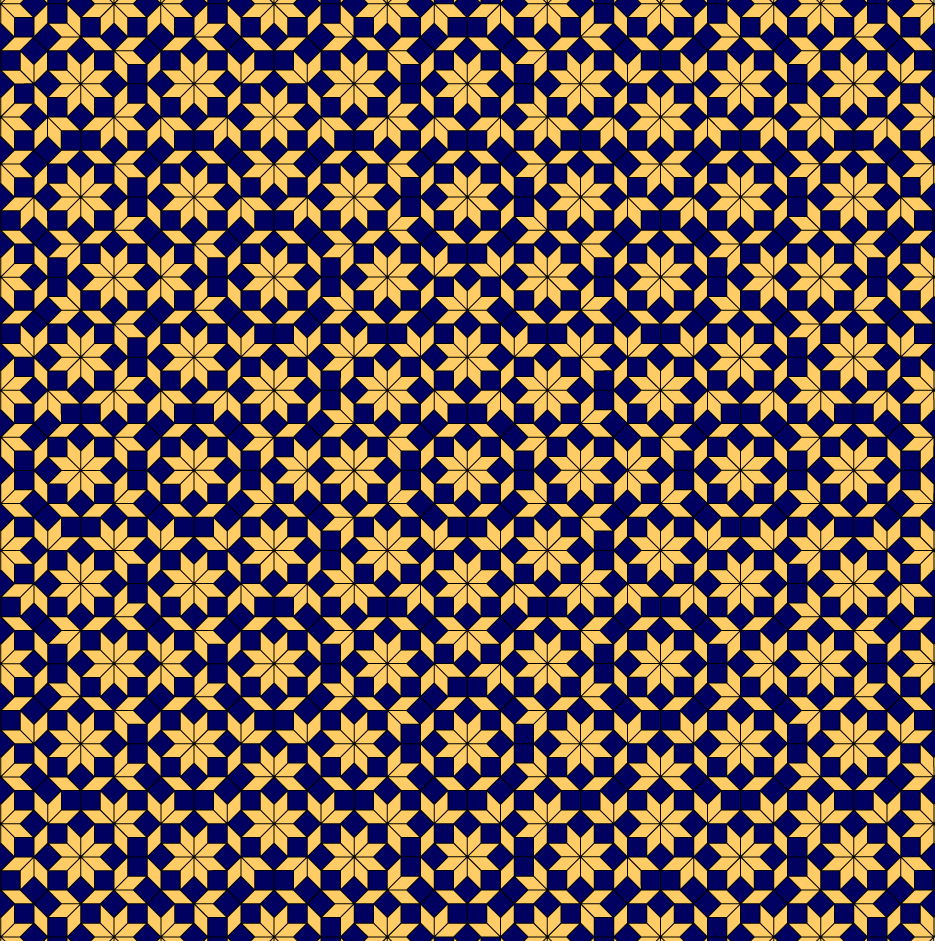}}
\caption{ The Ammann-Beenker or ``octahedral''  tiling.}
\label{FigAB1}
\end{figure}

Only last week (on May 4, 2018), R{\'e}my Sigrist investigated the
Ammann-Beenker (or ``octagonal'') tiling shown in Figure~\ref{FigAB1},
an aperiodic tiling with eight-fold rotational symmetry about the central point.

\begin{figure}[htb]
\centerline{\includegraphics[width=.60\textwidth]{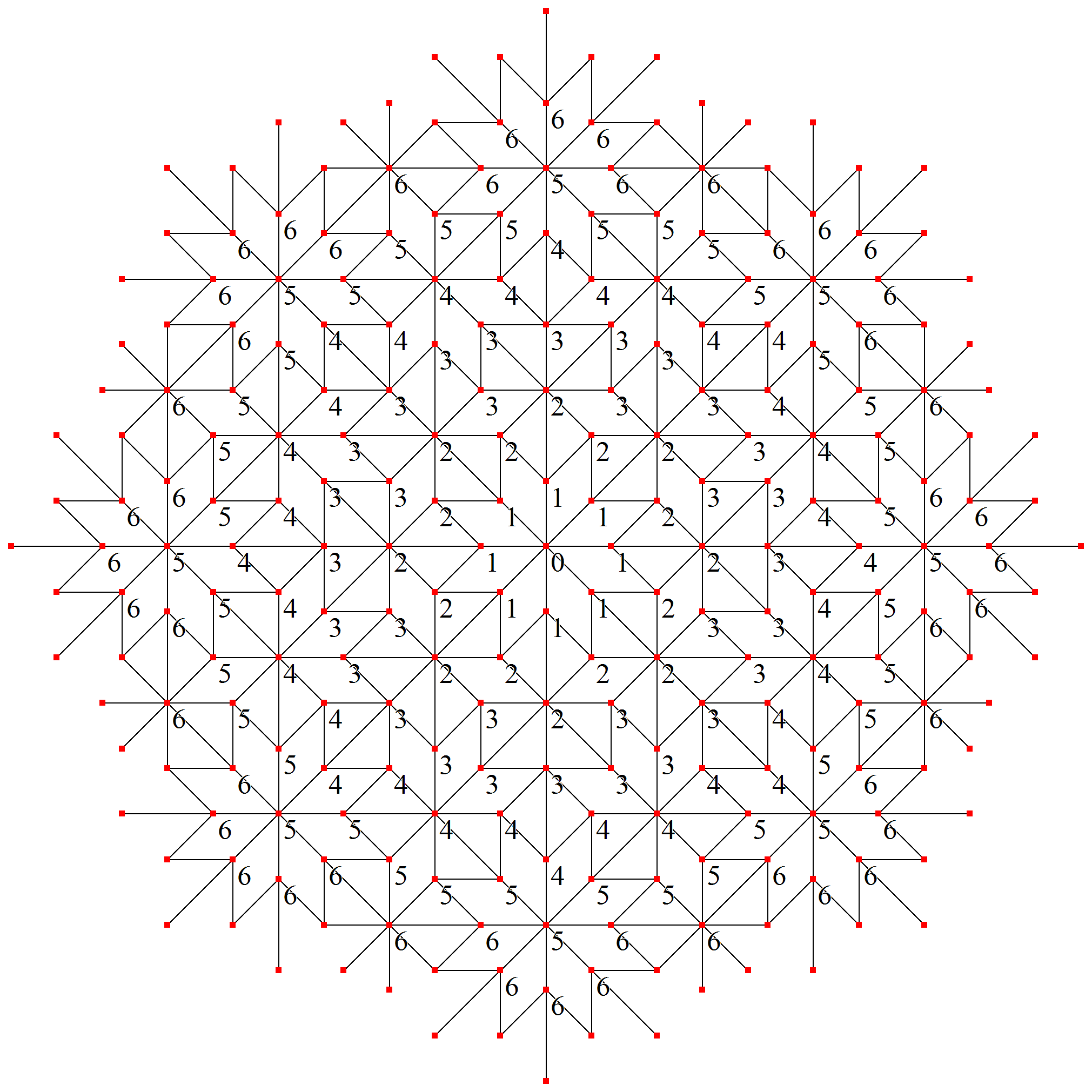}}
\caption{Illustrating the coordination sequence for the  Ammann-Beenker tiling, 
showing the vertices at distances $0$ though $6$ from the central vertex.}
\label{FigAB2}
\end{figure}

Sigrist determined the initial terms of the coordination sequence
with respect to the central vertex (\seqnum{A303981}):
\beql{EqAB1}
1, 8, 16, 32, 32, 40, 48, 72, 64, 96, 80, 104, 112, 112, 128, 152,  \ldots
\eeq
Figure~\ref{FigAB2} shows the vertices at distances $0,1,2,\ldots, 6$ from the center.

No formula or growth estimate is presently known for this sequence. However, earlier
this year  Anton Shutov and Andrey Maleev  determined the asymptotic behavior 
of the coordination sequence (\seqnum{A302176}) with respect to
a vertex with five-fold rotational symmetry in a certain Penrose tiling.
So we end with a question: Can the Shutov-Maleev approach
be used to find the asymptotic growth of \eqn{EqAB1}? 
Of course an explicit formula would be even nicer.

\section*{Acknowledgments}
Now that the OEIS is a wiki, many volunteer editors help maintain it.\footnote{And more are needed---if interested, please contact the author.}
Besides the OEIS contributors already mentioned in
this article, I would like to thank
J\"org Arndt,
Michael De Vlieger,
Charles Greathouse IV,
Alois P. Heinz,
Michel Marcus,
Richard J. Mathar,
Don Reble,
Jon E. Schoenfield,
and
Allan~C.~Wechsler
for their help.
(Many other names could be added to this list.)

\bigskip
\hrule
\bigskip


\noindent
ABOUT THE AUTHOR \\
Neil J. A. Sloane started collecting sequences in 1964. Every advance
in computing since then has made the task easier, and a subtitle could
be ``From punched cards to wiki in 54 years.'' 
The Unix operating system was the greatest breakthrough of all,
especially {\em grep}, the programmable editor {\em ex}, and pipes.

\bigskip
\hrule
\bigskip

\noindent
Image Credits 

\begin{itemize}
\item Figures~\ref{FigPQ5}-\ref{FigPQPK},
\ref{FigEKG1}, \ref{FigY1}, \ref{FigSpiral}, \ref{FigTag1}: Michael De Vlieger
(who also designed the ``Peace to the Max'' T-shirt in Figure~\ref{FigPQ11}).

\item
Figure \ref{FigWildJG}: Jessica Gonzalez.

\item
Figures \ref{FigWild4a}-\ref{FigWild5a}:   Jonathan Wild. 

\item
Figure \ref{FigY2}: Hans Havermann.

\item
Figure \ref{FigSpiralColor}:  Kerry Mitchell.

\item
Figure \ref{FigTag2}: Lars Blomberg.

\item
Figures \ref{FigCairo}, \ref{FigTet2}: Chaim Goodman-Strauss.

\item
Figure \ref{FigAB1}: The  Tilings Encyclopedia.

\item
Figure \ref{FigAB2}: R{\'e}my Sigrist.

\end{itemize}

\bigskip
\hrule
\bigskip

\end{document}